\documentclass[12pt]{amsart}

\usepackage{amssymb}  
\usepackage{latexsym} 
\usepackage{comment}

\newcommand{\Z}{{\mathbb Z}}
\newcommand{\Q}{{\mathbb Q}}
\newcommand{\F}{{\mathbb F}}
\newcommand{\BP}{{\mathbb P}}
\newcommand{\CO}{{\mathcal O}}

\newcommand{\cd}{\!\cdot\!}

\newenvironment{Proof}{\par\noindent{\sc Proof:}}%
                      {\hspace*{\fill}\nobreak$\Box$\par}

\newtheorem{Theorem}{Theorem}[section]
\newtheorem{Lemma}[Theorem]{Lemma}
\newtheorem{Proposition}[Theorem]{Proposition}

\numberwithin{equation}{section}

\begin{document}

\title[Rational squares at fixed distance from a fifth power]%
      {On the number of rational squares \\%
       at fixed distance from a fifth power}

\author{Michael Stoll}
\address{School of Engineering and Science,
         International University Bremen,
         P.O.Box 750561,
	 28725 Bremen, Germany.}
\email{m.stoll@iu-bremen.de}
\date{\today}


\begin{abstract}
  The main result of this note is that there are at most seven rational
  points (including the one at infinity) on the curve $C_A$ with the affine
  equation $y^2 = x^5 + A$ (where $A$ is a tenth power free integer)
  when the Mordell-Weil rank of the Jacobian of~$C_A$ is one.
  This bound is attained for $A = 18^2$.
\end{abstract}

\maketitle

\renewcommand{\baselinestretch}{1.1}
\renewcommand{\arraystretch}{1.3}


\section{Introduction}

Let $A \neq 0$ be a rational number. We are
interested in the number of rational solutions $(x,y)$ to the
equation $y^2 = x^5 + A$. In more geometric terms, this amounts to
counting the (affine) rational points on the curve~$C_A$ given by the
(affine) equation
\[ C_A : y^2 = x^5 + A \,. \]

In this note, we take up ideas from~\cite{Stoll} and apply them
to this family of genus~2 curves. In the following, $C_A$ will denote
a smooth projective model of the curve in question. With respect to
rational points, this means that there is one additional rational
point ``at infinity'', which we will denote~$\infty$ in what follows.

Let $J_A$ be the Jacobian of~$C_A$, and denote by~$r_A$ the Mordell-Weil
rank of~$J_A(\Q)$. Since $C_A$ and $C_B$ are isomorphic when the
quotient $A/B$ is a tenth power, we can (and will)
assume that $A$ is an integer, not divisible by the tenth
power of any prime. Let $n_A$ be half the number of ``non-trivial'' points
in~$C_A(\Q)$, i.e., finite points with non-vanishing $x$ and $y$ coordinates.
Then $\#C_A(\Q) = 2\,n_A + d_A$, where $d_A = 1$, $2$, $3$, or $4$ if
$A$ is neither a square nor a fifth power, a fifth power but $A \neq 1$,
a square but $A \neq 1$, or $A = 1$, respectively. Since $r_1 = 0$,
we have $d_A \le 3$ if $r_A \ge 1$.

The numbers we are interested in are $n_A$ and $\#C_A(\Q) = d_A + 2n_A$. 
The result we will prove is as follows.

\begin{Theorem} \label{Thm}
  Let $A$ be a tenth power free integer, and assume that $r_A = 1$.
  Then $n_A \le 2$ and therefore $\#C_A(\Q) \le 7$.
  Furthermore, $\#C_A(\Q) = 7$ only when $A = 18^2$.
\end{Theorem}

Let $S$ be the set of tenth power free integers. If we define 
\begin{align*}
  N(r) &= \max\{n_A : A \in S, r_A = r\} \,, \\
  B(r) &= \max\{\#C_A(\Q) : A \in S, r_A = r\} \,,
\end{align*}
then the theorem says that $N(1) = 2$, $B(1) = 7$. 

For $r = 0$, we obtain $N(0) = 1$, $B(0) = 4$. This is because the
torsion points on~$C_A$ are known to be (see for example~\cite{PoonenTors})
\[ \infty\,, \quad (-\zeta^k \sqrt[5]{A}, 0)\,, \quad
   (0, \pm\sqrt{A}) \quad \text{and} \quad
   (\zeta^k \sqrt[5]{4A}, \pm\sqrt{5A}) \,,
\]
where $\zeta$ is a primitive fifth root of unity. The only nontrivial
points in this list are of the form given last. But there is only
one value of~$A$, namely $A = 2^8 \cd 5^5$, such that this leads
to a pair of rational points on~$C_A$. We have $r_A = 0$ and $n_A = 1$
in this case (but $d_A = 1$, so $\#C_A(\Q) = 3$). For all other~$A$
such that $r_A = 0$, we must have $n_A = 0$, hence $\#C_A(\Q) = d_A \le 4$.
The maxmimum is attained for the unique~$A$ with $d_A = 4$, namely~$A = 1$.

Since the method of proof can be applied only when $r_A < g(C_A) = 2$,
we cannot obtain exact values for higher ranks. However, we have found
the following examples, thus obtaining lower bounds.

\[ \begin{array}{|c|c|c|l|l|} \hline
     r & N(r) & B(r) & \text{$A$ with max.~$N(r)$} 
                     & \text{$A$ with max.~$B(r)$} \\\hline
     0 & 1 & 4 & 2^8 \cd 5^5 & 1 \\
     1 & 2 & 7 & 2^2 \cd 3^4 & 2^2 \cd 3^4 \\
     2 & \ge 3 & \ge 9 & 2^2 \cd 3^4 \cd 7^4 & 2^2 \cd 3^4 \cd 7^4 \\
     3 & \ge 4 & \ge 11 & 2^2 \cd 3^2 \cd 5^4 \cd 7^4 
                        & 2^2 \cd 3^2 \cd 5^4 \cd 7^4 \\
     4 & \ge 6 & \ge 15 & 3^4 \cd 7^4 \cd 19^4 & 3^4 \cd 7^4 \cd 19^4 
     \\\hline
   \end{array}
\]

For these examples, the rank~$r_A$ was determined by first computing
an upper bound using $2$-descent as described in~\cite{StollDesc} and
then exhibiting sufficiently many independent points in the Mordell-Weil
group (which are here provided by the rational points on~$C_A$).


\section{The method}

We will apply Chabauty's method with a twist, as explained in~\cite{Stoll}.
In that paper, we were only considering sufficiently large primes. In our
situation, we obtain the following result in this way. In the following,
$v_p$ denotes the normalized $p$-adic valuation.

\pagebreak[2]

\begin{Proposition} \label{Prop}
  Suppose that $r_A = 1$. Then $n_A \le 1$ in the following cases.
  \begin{itemize}\addtolength{\itemsep}{1mm}
    \item[(1)] $v_p(A) = 5$ for some prime $p \ge 7$.
    \item[(2)] $v_p(A) \in \{2,4,6,8\}$ for some prime $p \ge 11$.
    \item[(3)] $v_p(A) \in \{1,3,7,9\}$ for some prime $p \ge 17$.
  \end{itemize}
\end{Proposition}
\begin{Proof}
  We use the main result of~\cite{Stoll}, applied to $C_A$ with
  $\Gamma$ taken to be (1) $\mu_2$ acting on~$y$, (2) $\mu_5$ acting on~$x$,
  and (3) $\mu_2 \times \mu_5$ acting on~$(y,x)$, respectively.
  Note that in each case, $C_A$ is a $\Gamma$-twist of $C_{A/p^{v_p(A)}}$,
  which has good reduction at~$p$.
\end{Proof}

In principle, this reduces the cases that we have to check to a finite number.
However, the number of cases is large (a priori, there are
$10^3 \cd 9 \cd 5^2 = 225\,000$ curves; we can expect close to half
of them to have $r_A = 1$), and dealing with them one by one would
require a very large amount of computation.

We therefore want to use the method at the small primes, too.
This will reduce the cases we have to look at to a manageable number.

The basic setup is as follows.
If $r_A = 1$, then there is a differential
\[ \omega = \frac{(\alpha + \beta\,x)\,dx}{2y} \in \Omega(C_A/\Q_p) \]
(with $(\alpha : \beta) \in \BP^1(\Q_p)$) killing the Mordell-Weil group
in the sense that
\[ \lambda_\omega(P) = \int\limits_0^P \omega = 0 
     \quad\text{for all $P \in J_A(\Q)$} \,. 
\]
Note that the integral is linear both in~$P \in J_A(\Q_p)$ and in~$\omega$, and
vanishes when $P$ is a torsion point in~$J_A(\Q_p)$.
We embed $C_A$ into~$J_A$ using the point at infinity as a base point
and from now on consider $C_A$ as a subvariety of~$J_A$.
Then $C_A(\Q)$ is contained in the set of zeros of~$\lambda_\omega$
on~$C_A(\Q_p)$. Therefore the number of nontrivial rational points
is bounded by the number of nontrivial zeros of~$\lambda_\omega$
(note that $\lambda_\omega$ vanishes at the trivial points, since they
are mapped to torsion points of~$J_A$).

Before we use the twisting trick of~\cite{Stoll}, let us prove a result that
will help reduce the number of cases later, using the standard Chabauty
technique.

\begin{Lemma} \label{Lem30}
  Assume that $r_A = 1$. If $A \equiv 1 \bmod 3$, then we have $n_A \le 1$. 
  If $A \equiv -1 \bmod 3$, then we have $n_A \le 2$.
\end{Lemma}

\begin{Proof}
  If $A \equiv 1 \bmod 3$, then we have 
  $C_A(\F_3) = \{\infty, (0, 1), (0, -1), (-1, 0)\}.$
  If $A \equiv -1 \bmod 3$, then
  $C_A(\F_3) = \{\infty, (-1, 1), (-1, -1), (1, 0)\}.$
  
  Let $\bar{\omega}$ be the reduction of (a suitable multiple of) $\omega$
  mod~$3$. We find that $v(3, 0) = 0$, $v(3, 1) = 1$, $v(3, 2) = 0$ in the
  notation of~\cite[\S\,6]{Stoll}. By Prop.~6.3 in~\cite{Stoll},
  the number of zeros of~$\lambda_\omega$ on the residue class 
  of~$P \in C_A(\F_3)$ is at most $1 + n + v(3,n)$, where
  $n = v_P(\bar{\omega})$. 
  This implies that we can only get nontrivial points in residue
  classes on which $\bar{\omega}$ vanishes, or in nontrivial residue classes.
  Furthermore, the number of zeros of $\lambda_\omega$ in the residue
  class of~$P$ can be at most~$3$, 
  since $v_P(\bar{\omega}) \le 2$.
  The two nontrivial classes $(-1, \pm 1)$ that occur for 
  $A \equiv -1 \bmod 3$ contain the torsion point
  $(\sqrt[5]{4A}, \pm\sqrt{5A}) \in C_A(\Q_3)$. Since this point is
  rational only when $A = 2^8 \cd 5^5$, and $r_A = 0$ in this case,
  it is still true that we can get nontrivial rational points in these
  classes only when $\bar{\omega}$ vanishes there.

  Let us consider the various classes in turn. 
  If $P = \infty$ or $P = (\pm 1, 0)$ and $\bar{\omega}$ vanishes at~$P$
  (then of second order), there can only be one pair of nontrivial rational
  points in this residue class, and hence $n_A \le 1$, since no other
  class contributes to~$n_A$.
  
  In the other cases, $\bar{\omega}$ can only vanish to first order. 
  One of the a priori up to three zeros on the residue class will be a
  torsion point, which is trivial or not rational. Therefore
  there can be at most two pairs of nontrivial rational points, and $n_A \le 2$.
  It remains to show that in fact, $n_A \le 1$ if $A \equiv 1 \bmod 3$.
  In this case, we have $P = (0, \pm 1)$ with $v_P(\bar{\omega}) = 1$,
  so $\bar{\omega} = x\,dx/2y$.
  
  We write $A = a^2$ with $a \in \Z_3^\times$ such that $(0, a)$ reduces 
  to the point we are
  considering. $t = x$ is a uniformizer, and we have, taking
  $\omega = (3\alpha + x)dx/2y$:
  \begin{align*}
    \frac{y}{a} &= 1 + \frac{t^5}{2A} + O(t^{10}) \\
    \frac{a\,dx}{y} &= \Bigl(1 - \frac{t^5}{2A} + O(t^{10})\Bigr)\,dt \\
    2a\,\omega
      &= \Bigl(3\alpha + t - \frac{3\alpha}{2A}\,t^5 + O(t^7)\Bigr)\,dt \\
    2a\,\lambda_\omega
      &= 3\alpha\,t + \frac{t^2}{2} - \frac{\alpha}{4A}\,t^6 + O(t^8)
  \end{align*}
  We see that there is only one non-trivial root in $3\Z_3$ (satisfying
  $t \equiv 3\alpha \bmod 9$); therefore $n_A \le 1$.
\end{Proof}


\section{Looking at the small primes}

Our first goal is to show that $n_A \le 2$ for all~$A$ such that $r_A = 1$,
by a detailed study of the $3$-adic situation. We have covered the
case $v_3(A) = 0$ in Lemma~\ref{Lem30} already, so here we will
consider the cases $1 \le v_3(A) \le 9$.

We keep the notations introduced above. However, we now suppose that
$A$ is divisible by~$p$ (we will mostly take $p = 3$ below), 
so that $A = p^\nu\,a$ with some
$1 \le \nu \le 9$ and $a \in \Z_p^\times$. 
Let $\pi = p^{1/10}$ 
and set $x = \pi^{2\nu}\,X$, $y = \pi^{5\nu}\,Y$;
then over $\Q_p(\pi)$, $C_A$ is isomorphic to
\[ C_a : Y^2 = X^5 + a \,, \]
and on~$C_a$,
\[ \omega = (\pi^{-\nu}\alpha + \pi^\nu \beta X) \frac{dX}{2Y} \,. \]
Since $\Q_p(\pi)/\Q_p$ is totally ramified, the residue class field 
of~$\Q_p(\pi)$ is~$\F_p$.
The points in $C_A(\Q_p)$ are mapped to one of the following 
types of points 
in~$C_a(\F_p)$. ($\bar{b} \in \F_p$ denotes the image of~$b \in \Z_p$.)
\begin{itemize}\addtolength{\itemsep}{1mm}
  \item[1.] $\infty$.
  \item[2.] $(-\bar{b}, 0)$ if $\nu = 5$ and $a = b^5$ for some $b \in \Z_p$.
  \item[3.] $(0, \bar{b})$ if $\nu$ is even and $a = b^2$ for some $b \in \Z_p$.
\end{itemize}
This holds when $p \ne 2, 5$. It still holds when $p = 2$ or~$5$ and
$p \nmid \nu$.

We want to bound the number of non-trivial points in~$C_A(\Q)$ mapping
to each of these points mod~$\pi$. If $p$ is large enough, this bound is
given by the order of vanishing of the differential $\bar{\omega}$ at
the point in question (assuming $\omega$ is scaled such that it is
integral and reduces to something non-zero mod~$\pi$). This is how
the results in Prop.~\ref{Prop} are obtained. In order to get bounds
when $p$ is small, we need to take a closer look at the logarithm
\[ \lambda_\omega(T) := \lambda_\omega(P(T))
      = \int\limits_{0}^{P(T)} \omega
      = \int\limits_{P(0)}^{P(T)} \omega \,, \]
(recall that the logarithm vanishes on torsion points)
where $T$ is a uniformizer at the trivial point $P(0)$ in the residue
class under consideration, and $P(T)$ is the point corresponding to 
the value $T \in \pi \CO_\pi$ of the uniformizer. This logarithm
$\lambda_\omega$ can be expanded into a power series in~$T$, and
the number of its zeros in $\pi \CO_\pi$ can be bounded above by
considering the valuations of the coefficients (and, in some cases, 
the factoring of polynomials over~$\F_p$). In fact, we are only interested
in zeros that arise from points in~$C_A(\Q_p)$, which restricts the
possibilities further.

We discuss the various possible image points in turn.

\subsection*{The residue class of $\infty$} \strut \\
The corresponding points on~$C_A$ have $t = x^2/y \in \Z_p$, and we have
$T = X^2/Y = \pi^\nu\,t \in \pi^\nu\,\Z_p$. $T$ is a uniformiser in~$\infty$
on~$C_a$. We have the equations
\[ X^{-1} = T^2\,(1 + a\,X^{-5})\,, \qquad Y^{-1} = T\,X^{-2} \,, \]
so
\begin{align*}
   X      &= T^{-2}(1 - a\,T^{10} + O(T^{20})) \\
   -\frac{dX}{2Y} &= T^2\,(1 + 6a\,T^{10} + O(T^{20}))\,dT \\
   -\omega
          &= (\beta\,\pi^{\nu} + \alpha\,\pi^{-\nu}\,T^2 
                + 5a\beta\,\pi^{\nu}\,T^{10}
                + 6a\alpha\,\pi^{-\nu}\,T^{12} + O(T^{20}))\,dT \\
   -\lambda_\omega
          &= \beta\,\pi^{\nu}\,T + \frac{\alpha}{3}\,\pi^{-\nu}\,T^3
              + \frac{5a\beta}{11}\pi^{\nu}\,T^{11}
              + \frac{6a\alpha}{13}\pi^{-\nu}\,T^{13} + O(T^{21}) \\
          &= \pi^{2\nu}\Bigl(\beta\,t + \frac{\alpha}{3}\,t^3
                             + \frac{5a\beta}{11}\,p^\nu\,t^{11}
                             + \frac{6a\alpha}{13}\,p^\nu\,t^{13} + \dots\Bigr)
\end{align*}
Since $3$, $11$ and~$13$ are the only primes occurring in denominators
of relevant coefficients
(the later terms do not matter, as is easily seen), we see that for all
other primes~$p$, the
following holds.

If $\overline{(\alpha : \beta)} = (0 : 1) \in \BP^1(\F_p)$, then there
is only one rational point (namely~$\infty$) in this residue class. Otherwise,
if $\overline{(\alpha : \beta)} = (1 : 0)$, there may be three, and if
$\overline{(\alpha : \beta)} = (1 : \xi)$ with $\xi \neq 0$, there is
one point if $-3\alpha\beta$ is a non-square mod~$p$, and at most 
three points if $-3\alpha\beta$ is a square mod~$p$.

We will not discuss $p = 11$ and $p = 13$ here.
If $p = 3$, there always are at most three points, but the condition is
shifted. We write $\omega = (3\alpha' + \beta\,x)\,dx/2y$ with 
$\alpha', \beta \in \Z_3$; then there can
be three points if $\bar{\alpha}' \neq 0$ and $-\bar{\alpha}'\bar{\beta}$ is
a square. We obtain the following result, strengthening part~(3) of
Prop.~\ref{Prop}.

\begin{Lemma} \label{Lem31}
  Suppose $r_A = 1$. If $v_p(A) \in \{1,3,7,9\}$ for some $p \neq 11, 13$, 
  then $n_A \le 1$.
\end{Lemma}

\begin{Proof}
  Note that if $v_p(A) \in \{1,3,7,9\}$, $\infty$ is the only point 
  in~$C_a(\F_p)$ that is hit by~$C_A(\Q)$. By the preceding discussion,
  there is at most one pair of nontrivial rational points in this
  residue class.
\end{Proof}

\subsection*{The residue class of $(-b, 0)$} \strut \\
we now assume $\nu = 5$ (and $p \neq 5$)
and $a = b^5$, so we have $y^2 = x^5 + p^5 b^5$.
The points in the residue class we are considering have 
$x \equiv -bp \bmod p^2$ and $y \equiv 0 \bmod p^3$. We choose $T = Y$
as the uniformiser on~$C_a$; then $T = \sqrt{p}\,t$ with $t \in \Z_p$.
Expanding everything in terms of~$T$, we get
\begin{align*}
  X &= -b(1 - a^{-1}\,T^2)^{1/5} \\
    &= -b\Bigl(1 - \frac{1}{5a}\,T^2 - \frac{2}{(5a)^2}\,T^4
                 - \frac{6}{(5a)^3}\,T^6 
                 + O(T^8)\Bigr) \\
  \frac{5a}{b}\,\frac{dX}{2Y}
    &= \Bigl(1 + \frac{4}{5a}\,T^2 + \frac{18}{(5a)^2}\,T^4
               + O(T^8)\Bigr)\,dT \\
  \frac{5a}{b}\,\frac{X\,dX}{2Y}
    &= -b\Bigl(1 +\frac{3}{5a}\,T^2 + \frac{12}{(5a)^2}\,T^4
                 + O(T^8)\Bigr)\,dT \\
  \frac{5a}{b}\,\omega
    &= \frac{1}{\sqrt{p}}\Bigl((\alpha - b \beta p)
                 + \frac{4 \alpha - 3 b \beta p}{5a}\,T^2
                 + \dots\Bigr)\,dT \\
  \frac{5a}{b}\,\lambda_\omega
    &= \frac{1}{\sqrt{p}}(\alpha - b \beta p)\,T 
                 + \frac{4 \alpha - 3 b \beta p}{3 \cd 5a}\,T^3
                 + \dots \\
    &= (\alpha - b \beta p)\,t
                      + \frac{(4 \alpha - 3 b \beta p)p}{3 \cd 5a}\,t^3
                      + \dots
\end{align*}
The interesting case for us here is $p = 3$. 
We again write $\alpha = 3\alpha'$ and assume that $\alpha', \beta \in \Z_3$,
not both in~$3\Z_3$; 
then up to scaling, we have mod~$3$
\[ \lambda_\omega \sim t\,((\bar{\alpha}' - \bar{b}\,\bar{\beta}) 
                        - \bar{a}^{-1}\,\bar{\alpha}'\,t^2) \,. 
\]
If $\bar\alpha' = 0$, there will be only one solution. If $\alpha' = 1$,
we will have three solutions if $a(1 - b \beta)$ is a square, and one
solution if it is a non-square. In any case, there is at most one pair
of nontrivial rational points in this residue class. 
For all other primes $p \neq 5$, we get at most one pair of nontrivial
rational points in this class as well, but only if $v_p(\alpha) > v_p(\beta)$.
We therefore obtain,
taking into account the discussion of~$\infty$ above, the following 
strengthening of Prop.~\ref{Prop}, part~(1).

\begin{Lemma} \label{Lem35}
  Suppose $r_A = 1$. If $v_3(A) = 5$, then $n_A \le 2$.
  If $v_p(A) = 5$ for a prime $p \notin \{3, 5\}$, then $n_A \le 1$.
\end{Lemma}

\subsection*{The residue class of $(0, b)$} \strut \\
Let $\nu = 2n$ and set $\mu = \min\{ m \in \Z \mid 5m > n \}$. 
We must exclude $p = 2$ here. We have
$\rho = 5\mu - \nu = 3, 1, 4, 2$ for $n = 1, 2, 3, 4$, respectively.
The points in the residue class have $x = p^\mu\,t$ with $t \in \Z_p$ and 
$y \equiv b p^n \bmod p^{n+1}$. We choose $T = X = \pi^{2\rho}\,t$ 
as a uniformiser on~$C_a$.
This gives
\begin{align*}
  Y^{-1} &= b^{-1}\,(1 + a^{-1}\,T^5)^{-1/2} 
          = b^{-1}\Bigl(1 - \frac{1}{2a}\,T^5 + \frac{3}{8a^2}\,T^{10}
                          + O(T^{15})\Bigr) \\
  b\,\frac{dX}{Y}
         &= \Bigl(1 - \frac{1}{2a}\,T^5 + \frac{3}{8a^2}\,T^{10}
                    + O(T^{15})\Bigr)\,dT \\
  2b\,\omega
         &= \Bigl(\alpha\,\pi^{-2n} + \beta\,\pi^{2n}\,T
                   - \frac{\alpha}{2a}\pi^{-2n}\,T^5 
                   - \frac{\beta}{2a}\pi^{2n}\,T^6
                   + O(T^{10})\Bigr)\,dT \\
  2b\,\lambda_\omega
         &= \alpha\,\pi^{-2n}\,T + \frac{\beta}{2}\pi^{2n}\,T^2
              - \frac{\alpha}{3 \cd 4a}\pi^{-2n}\,T^6
              - \frac{\beta}{7 \cd 2a}\pi^{2n}\,T^7 + \dots \\
         &= \pi^{2(\rho-n)} \Bigl(\alpha\,t + \frac{\beta}{2}\,p^\mu\,t^2
                               - \frac{\alpha}{3 \cd 4a}\,p^\rho\,t^6
                               - \frac{\beta}{7 \cd 2a}\,p^{\rho+\mu}\,t^7
                               + \dots \Bigr)
\end{align*}

If $p \notin \{2,3,7\}$, then there is at most one nontrivial solution,
and we obtain at most one pair of nontrivial rational points mapping
to $(0, \pm\bar{b})$. Together with the discussion of~$\infty$,
this proves the following. The only new case is $p = 5$,
since $p \ge 11$ is already taken care of by Prop.~\ref{Prop}, part~(2).

\begin{Lemma} \label{Lem5}
  Assume $r_A = 1$ and $v_p(A) \in \{2,4,6,8\}$ for some $p \notin \{2,3,7\}$.
  Then $n_A \le 1$.
\end{Lemma}

\begin{Proof}
  We have a contribution of at most~$1$ to~$n_A$ from~$\infty$  and a 
  contribution of at most~$1$ from $(0, \pm \bar{b})$. However, we get
  a contribution from~$\infty$ only when $v_p(\beta) \ge v_p(\alpha)$,
  but then there is no contribution from $(0, \pm \bar{b})$, as can be
  seen from the expansion of $\lambda_\omega$ above.
\end{Proof}

We have to consider the cases $p = 3$ and $p = 7$ separately. 
When $p = 3$, we have, with $\omega = (3\alpha' + \beta x)dx/2y$ as before:
\[ \lambda_\omega \sim \alpha'\,t + \frac{\beta}{2}\,3^{\mu-1}\,t^2
             - \frac{\alpha'}{4a}\,3^{\rho-1}\,t^6 + \dots \,. \]
If $\rho > 1$, only the first two terms matter. We get extra
solutions only if $v_3(\alpha') \ge v_3(\beta) + \mu - 1$.
If $\rho = 1$ (and therefore $\mu = 1$, $\nu = 4$), we have to look
at solutions in~$\F_3$ of
\[ t(\bar\alpha' - \bar{\beta}\,t - \bar\alpha'\,t^5) \,. \]
(Recall that $a$ is a square, so $\bar{a} = 1$.)
If $\bar\alpha' = 0$, there is one extra solution. Otherwise, we can
take $\alpha' = 1$, and then we have no extra solutions if
$\bar{\beta} = -1$, we have one extra solution
if $\bar{\beta} = 0$, and we have potentially two extra solutions
if $\bar{\beta} = 1$.
Considering this together with the result at infinity, we get the
following.

\begin{Lemma} \label{Lem32}
  Suppose $r_A = 1$. If $v_3(A) \in \{6, 8\}$, then $n_A \le 1$.
  If $v_3(A) \in \{2, 4\}$, then $n_A \le 2$.
\end{Lemma}

\begin{Proof}
  The following table summarizes the possible contributions to~$n_A$
  from the residue classes of~$\infty$ and of~$(0, \pm b)$, depending
  on the reduction mod~$3$ of $(\alpha' : \beta) \in \BP^1(\Q_3)$.
  \[ \begin{array}{|c||c|ccc|} \hline
       \overline{(\alpha' : \beta)} & \infty & \multicolumn{3}{|c|}{(0,\pm b)}
        \\
        & & \nu = 2 & \nu = 4 & \nu = 6,8 \\\hline
       (0 : 1)  &     0 & \le 1 & \le 1 & \le 1 \\
       (1 : 0)  & \le 1 &     0 & \le 1 &     0 \\
       (1 : 1)  &     0 & \le 1 & \le 2 &     0 \\
       (1 : -1) & \le 1 & \le 1 &     0 &     0 \\\hline
     \end{array}
  \]
  We see that for $\nu = 6$ or~$8$, $n_A \le 1$, whereas for $\nu = 2$ or~$4$,
  $n_A \le 2$.
\end{Proof}

Now we consider $p = 7$. This leads to
\[ \lambda_\omega \sim t(\alpha + \frac{\beta}{2}\,7^\mu\,t 
             - \frac{\alpha}{12a}\,7^\rho\,t^5
             - \frac{\beta}{2a}\,7^{\rho+\mu-1}\,t^6 + \dots) \,. \]
The $t^5$ term is irrelevant. If $\bar\alpha \neq 0$, then there are
no extra solutions. In general, we need $v_7(\alpha) \ge v_7(\beta) + \mu$
for there to be extra solutions. In this case, if $\rho > 1$, then there
is just one extra solution. If $\rho = 1$, i.e., $\nu = 4$ and 
$\alpha = 7 \alpha'$, then (taking $\beta = 1$ without loss of generality),
we must consider the roots in~$\F_7$ of
\[ \bar{\alpha}'\,t - 3\,t^2 + 3\bar{a}^{-1}\,t^7 \,. \]
There can be more than one extra solution; in this case there are up to four
extra solutions. In any case, together with the result
at infinity, we get the following.

\begin{Lemma} \label{Lem7}
  Suppose $r_A = 1$. If $v_7(A) \in \{2, 6, 8\}$, then $n_A \le 1$.
\end{Lemma}

\subsection*{Putting it Together}\strut

Collecting the information obtained so far, we see that Lemmas \ref{Lem30},
\ref{Lem31}, \ref{Lem35} and~\ref{Lem32} cover all cases. This proves
the first part of Theorem~\ref{Thm}.

Now, if $A$ is such that $r_A = 1$ and $\#C_A(\Q) = 7$, we need to have
$d_A = 3$ and $n_A = 2$. So $A$ has to be a square. Furthermore, 
by Lemma~\ref{Lem5}, the prime factors of~$A$ are contained in~$\{2,3,7\}$, 
by Lemma~\ref{Lem7},
$v_7(A) \in \{0, 4\}$, and by Lemmas \ref{Lem30} and~\ref{Lem32},
$v_3(A) \in \{2, 4\}$. This leaves $5 \cd 2 \cd 2 = 20$ values
of~$A$ to check. We can reduce the number of cases further by noting
that if $A \equiv 1, 3, 9 \bmod 11$, all points in~$C_A(\F_{11})$
lift to torsion points in~$C_A(\Q_{11})$, and therefore by standard
Chabauty, $n_A \le 1$. This reduces the list to the following eight values:
\[ A \in \{3^2 \cd 7^4,\; 3^4,\; 2^2 \cd 3^4,\; 2^2 \cd 3^4 \cd 7^4,\;
           2^4 \cd 3^4 \cd 7^4,\; 2^6 \cd 3^2, 2^8 \cd 3^2,\;
           2^8 \cd 3^2 \cd 7^4\}
\]
The table below summarizes the data for these curves.

\[ \begin{array}{|l|c|c||l|c|c|} \hline
      A & r_A & n_A & A & r_A & n_A \\\hline
      3^2 \cd 7^4         & 2 & \ge 1 &
      3^4                 & 2 & \ge 2 \\
      2^2 \cd 3^4         & 1 &     2 &
      2^2 \cd 3^4 \cd 7^4 & 2 & \ge 3 \\
      2^4 \cd 3^4 \cd 7^4 & 3 & \ge 2 &
      2^6 \cd 3^2         & 1 &     1 \\
      2^8 \cd 3^2         & 1 &     0 &
      2^8 \cd 3^2 \cd 7^4 & 0 &     0 \\
      \hline
  \end{array}
\]

\medskip

The ranks $r_A$ have been found by computing an upper bound as described
in~\cite{StollDesc} and exhibiting sufficiently many independent points
in~$J_A(\Q)$ (which are mostly provided by points in~$C_A(\Q)$). The values
for~$n_A$ in the cases when $r_A = 1$ have been verified by a standard
Chabauty computation (at $p = 29$ for $A = 2^8 \cd 3^2$, at $p = 29$
and $p = 59$ for $A = 2^6 \cdot 3^2$: there are four extra residue classes
left after the computation with $p = 29$, which can then be excluded
by looking mod~$59$ --- see~\cite[\S\,12]{PSS} for an explanation of the
method).

This shows that $A = 18^2 = 2^2 \cd 3^4$ is the only value such that
$r_A = 1$ and $C_A$ has $7$ rational points. The last statement
of Theorem~\ref{Thm} is therefore also verified.

It would be interesting to find out if there are more values of~$A$ such
that $r_A = 1$ and $n_A = 2$ (which then will be non-squares).

In any case, it is easy to see that there is no~$A$ such that $r_A = 1$
and $\#C_A(\Q) = 6$. This would imply that $n_A = 2$ and $A$ is a fifth
power, so by Lemma~\ref{Lem35}, $A$ is one of $1, 3^5, 5^5, 3^5 \cd 5^5$.
But all these values satisfy $A \equiv 1 \bmod 11$, so $n_A \le 1$ if
$r_A = 1$ (which is the case for $A = 3^5$ and $3^5 \cd 5^5$).


\end{document}